\documentclass[12pt]{amsart}
\usepackage{latexsym}
\title[$\mathbf{P}^2$-irreducible open 3-manifolds]
{$\mathbf{R}^2$-irreducible universal covering spaces of  
$\mathbf{P}^2$-irreducible open 3-manifolds}
\author{Robert Myers} 
\address{Department of Mathematics, Oklahoma State University, Stillwater, OK 74078}
\email{myersr@math.okstate.edu}
\thanks{Research at MSRI is supported in part by NSF grant DMS-9022140.}
\subjclass{Primary: 57M10; Secondary: 57N10, 57M60}
\keywords{3-manifold, open 3-manifold, $\mathbf{P}^2$-irreducible, 
$\mathbf{R}^2$-irreducible, universal covering space, contractible open 
3-manifold, Whitehead manifold}

\newtheorem{thm}{Theorem}
\newtheorem{conj}{Conjecture}
\newtheorem{lem}{Lemma}[section]
\newcommand{\R}{\ensuremath{\mathbf{R}}}
\newcommand{\RR}{\ensuremath{\mathbf{R}^2}}
\newcommand{\RRR}{\ensuremath{\mathbf{R}^3}}
\newcommand{\PP}{\ensuremath{\mathbf{P}^2}}
\newcommand{\bd}{\ensuremath{\partial}}
\newcommand{\irr}{irreducible}
\newcommand{\rirr}{\RR-\irr}
\newcommand{\pirr}{\PP-\irr}
\newcommand{\birr}{\bd-\irr}
\newcommand{\inc}{incompressible}
\newcommand{\tm}{3-manifold}
\newcommand{\Mt}{\ensuremath{\widetilde{M}}}
\newcommand{\Qt}{\ensuremath{\widetilde{Q}}}
\newcommand{\hm}{homeomorphic}
\newcommand{\8}{\ensuremath{\infty}}
\newcommand{\n}{^{-1}}
\newcommand{\p}{^{\prime}}
\newcommand{\inte}{int \, }

\newcommand{\halfspace}{\ensuremath{\RR \times [0,\8)}}
\newcommand{\prodspace}{\ensuremath{\RR \times [0,1]}}

\newcommand{\Dn}{\ensuremath{D_n}}
\newcommand{\Dnn}{\ensuremath{D_{n+1}}}
\newcommand{\Pn}{\ensuremath{P_n}}
\newcommand{\Qn}{\ensuremath{Q_n}}
\newcommand{\Rn}{\ensuremath{R_n}}
\newcommand{\Hn}{\ensuremath{H_n}}
\newcommand{\Vn}{\ensuremath{V_n}}
\newcommand{\Mn}{\ensuremath{M_n}}
\newcommand{\Pnn}{\ensuremath{P_{n+1}}}
\newcommand{\Qnn}{\ensuremath{Q_{n+1}}}
\newcommand{\Rnn}{\ensuremath{R_{n+1}}}
\newcommand{\Hnn}{\ensuremath{H_{n+1}}}

\newcommand{\Mnn}{\ensuremath{M_{n+1}}}
\newcommand{\Pnj}{\ensuremath{P_{n,j}}}
\newcommand{\Rnj}{\ensuremath{R_{n,j}}}
\newcommand{\Hnj}{\ensuremath{H_{n,j}}}
\newcommand{\Vnj}{\ensuremath{V_{n,j}}}
\newcommand{\Simn}{\ensuremath{\Sigma^m_n}}
\newcommand{\Lamn}{\ensuremath{\Lambda^m_n}}
\newcommand{\Phmn}{\ensuremath{\Phi^m_n}}
\newcommand{\Cmn}{\ensuremath{C^m_n}}
\newcommand{\Cmnn}{\ensuremath{C^m_{n+1}}}
\newcommand{\Pnnj}{\ensuremath{P_{n+1,j}}}

\newcommand{\Rnnj}{\ensuremath{R_{n+1,j}}}
\newcommand{\Hnnj}{\ensuremath{H_{n+1,j}}}

\newcommand{\alz}{\ensuremath{\alpha_0}}
\newcommand{\alo}{\ensuremath{\alpha_1}}
\newcommand{\alt}{\ensuremath{\alpha_2}}
\newcommand{\beo}{\ensuremath{\beta_1}}
\newcommand{\bet}{\ensuremath{\beta_2}}
\newcommand{\gao}{\ensuremath{\gamma_1}}
\newcommand{\gat}{\ensuremath{\gamma_2}}
\newcommand{\deo}{\ensuremath{\delta_1}}
\newcommand{\detw}{\ensuremath{\delta_2}}
\newcommand{\ep}{\ensuremath{\varepsilon}}
\newcommand{\tht}{\ensuremath{\theta}}
\newcommand{\tho}{\ensuremath{\theta_1}}
\newcommand{\thz}{\ensuremath{\theta_0}}
\newcommand{\om}{\ensuremath{\omega}}
\newcommand{\omj}{\ensuremath{\omega_j}}
\newcommand{\pp}{^{\prime\prime}}
\newcommand{\lao}{\ensuremath{\lambda_1}}
\newcommand{\lat}{\ensuremath{\lambda_2}}
\newcommand{\kao}{\ensuremath{\kappa_1}}
\newcommand{\kat}{\ensuremath{\kappa_2}}
\newcommand{\kad}{\ensuremath{\kappa_3}}
\newcommand{\nul}{\emptyset}
\newcommand{\pe}{poly-excellent}
\newcommand{\Pnnjp}{\ensuremath{P_{n+1,j}^+}}
\newcommand{\Simnn}{\ensuremath{\Sigma^m_{n+1}}}
\newcommand{\Lamnn}{\ensuremath{\Lambda^m_{n+1}}}
\newcommand{\Omn}{\ensuremath{O^m_n}}
\newcommand{\Icm}{\ensuremath{I_m}}
\newcommand{\Fmp}{\ensuremath{F_m\p}}
\newcommand{\Fmpp}{\ensuremath{F_m\pp}}
\newcommand{\Fmm}{\ensuremath{F_{m+1}}}
\newcommand{\Fm}{\ensuremath{F_m}}
\newcommand{\Fmmp}{\ensuremath{F_{m+1}\p}}
\newcommand{\Ft}{\ensuremath{\widetilde{F}}}
\newcommand{\Ej}{\ensuremath{E_j}}
\newcommand{\ali}{\ensuremath{\alpha_i}}
\newcommand{\alk}{\ensuremath{\alpha_k}}
\newcommand{\bei}{\ensuremath{\beta_i}}
\newcommand{\bek}{\ensuremath{\beta_k}}
\newcommand{\gai}{\ensuremath{\gamma_i}}
\newcommand{\gak}{\ensuremath{\gamma_k}}
\newcommand{\dei}{\ensuremath{\delta_i}}
\newcommand{\dek}{\ensuremath{\delta_k}}
\newcommand{\Pin}{\ensuremath{P^i_n}}
\newcommand{\Qin}{\ensuremath{Q^i_n}}
\newcommand{\Hin}{\ensuremath{H^i_n}}
\newcommand{\lait}{\ensuremath{\lambda^i_t}}
\newcommand{\Limun}{\ensuremath{L^{i,\mu}_n}}
\newcommand{\Zimun}{\ensuremath{Z^{i,\mu}_n}}
\newcommand{\Zimu}{\ensuremath{Z^{i,\mu}}}
\newcommand{\Cimum}{\ensuremath{C^{i,\mu,m}}}
\newcommand{\Cimumn}{\ensuremath{C^{i,\mu,m}_n}}
\newcommand{\Cmmn}{\ensuremath{C^{m+1}_n}}

\newcommand{\Mnnn}{\ensuremath{M_{n-1}}}
\newcommand{\Rnnn}{\ensuremath{R_{n-1}}}
\newcommand{\Cmnnn}{\ensuremath{C^m_{n-1}}}
\newcommand{\Xns}{\ensuremath{X_{n,s}}}
\newcommand{\Ks}{\ensuremath{K_{\sigma}}}
\newcommand{\Mtphi}{\ensuremath{\widetilde{M}[\varphi]}}
\newcommand{\Mtpsi}{\ensuremath{\widetilde{M}[\psi]}}
\newcommand{\Ut}{\ensuremath{\widetilde{U}}}

\begin{document}

\begin{abstract} An \irr\ open \tm\ $W$ is \rirr\ if it contains no  
non-trivial planes, i.e. given any proper embedded plane $\Pi$ in $W$ 
some component of $W-\Pi$ must have closure an  
embedded halfspace \halfspace. In this paper it 
is shown that if $M$ is a connected, \pirr, open \tm\ such that $\pi_1(M)$ 
is finitely generated and the universal covering space \Mt\ of $M$ is 
\rirr, then either \Mt\ is \hm\ to \RRR\ or $\pi_1(M)$ is a free product 
of infinite cyclic groups and fundamental groups of closed, connected 
surfaces other than $S^2$ or \PP. Given any finitely 
generated group $G$ of this form, uncountably many \pirr, 
open 3-manifolds $M$ are constructed with $\pi_1(M)\cong G$ such that the 
universal covering space \Mt\ is \rirr\ and not \hm\ to \RRR; the \Mt\ are 
pairwise non-\hm. Relations are 
established between these 
results and the conjecture that the universal covering space of any 
\irr, orientable, closed \tm\ with infinite fundamental group must be 
\hm\ to \RRR. 
\end{abstract}

\maketitle

\section{Introduction}

Suppose $M$ is a connected, \pirr, open \tm\ with $\pi_1(M)$ finitely 
generated and non-trivial. It is easy to construct examples of such $M$ 
for which the 
universal covering space \Mt\ is not \hm\ to \RRR. Start with any \tm\ 
$N$ satisfying the given conditions. Let $U$ be a \textbf{Whitehead 
manifold}, i.e. an \irr, contractible, open \tm\ which is not \hm\ 
to \RRR (see e.g. \cite{Wh}, \cite{Mc}). Choose end-proper embeddings of 
$[0,\8)$ in 
each of $N$ and $U$. (A map between manifolds is \textbf{end-proper} 
if pre-images of compact sets are compact; it is \textbf{$\bd$-proper} if the 
pre-image of the boundary is the boundary; it is \textbf{proper} if it has 
both these properties. These terms are applied to a submanifold if its 
inclusion map has the corresponding property.) Let $X$ and $Y$ be the 
exteriors of these rays. 
(The \textbf{exterior} of a submanifold is the closure of the complement of 
a regular neighborhood of it.) $\bd X$ and $\bd Y$ are each planes. We 
identify them to obtain a \pirr\ open \tm\ $M$ with $\pi_1(M) \cong 
\pi_1(N)$. Let $p:\Mt \rightarrow M$ be the universal covering map. 
Then \Mt, $p\n (X)$, and $p\n (Y)$ are \pirr\ \cite{MSY}. Each component 
$\widetilde{Y}$ of $p\n(Y)$ has interior \Ut\ \hm\ to $U$ and so contains a 
compact, connected subset $J$ which does not lie in a 3-ball in 
\Ut. If \Mt\ were \hm\ to \RRR\ then $J$ would lie in a 3-ball 
$B$ in \Mt. Standard general position and minimality arguments applied to 
$\bd B$ and $\bd \widetilde{Y}$ would then yield a 3-ball $B\p$ in \Ut\ 
containing $J$, a contradiction. Alternatively, one could use the Tucker 
Compactification Theorem \cite{Tu} to obtain a compact polyhedron $K$ in 
\Ut\ such that some component $V$ of $\Ut-K$ has 
non-finitely generated fundamental group. But this is impossible since 
the union of $V$ and $\Mt-\Ut$  
is a component of $\Mt-K$ whose fundamental group is isomorphic 
to $\pi_1(V)$. 

In this example $\bd \widetilde{Y}$ is a \textbf{non-trivial plane} 
in \Mt, i.e. a proper plane $\Pi$ such that no component of 
$\Mt-\Pi$ has closure \hm\ to \halfspace\ with $\Pi=\RR \times \{0\}$. 
This paper shows that it is harder to find examples if one rules out this 
behavior by requiring that \Mt\ be \textbf{\rirr} in the sense that, in 
addition to being \irr, it contains no non-trivial planes. 

Define a \textbf{closed surface group} to be the fundamental group of a 
closed, connected 2-manifold. 

\begin{thm} Let $M$ be a connected, \pirr, open \tm\ with $\pi_1(M)$ 
finitely generated. If the universal covering space \Mt\ of $M$ is 
\rirr, then either 
\begin{enumerate}
\item \Mt\ is \hm\ to \RRR\ or 
\item $\pi_1(M)$ is a free product of infinite cyclic groups and 
infinite closed surface groups. \end{enumerate} \end{thm}

The second possibility can be disjoint from the first. 

\begin{thm} Suppose $G$ is a free product of finitely many infinite 
cyclic groups and infinite closed surface groups. Then there is a 
\pirr\ open \tm\ $M$ such that $\pi_1(M)\cong G$ and $\Mt$ is an 
\rirr\ Whitehead manifold. Moreover, for each given $G$ there are 
uncountably many such $M$ for which the \Mt\ are pairwise non-\hm. 
\end{thm}

This generalizes an example of Scott and Tucker \cite{Sc-Tu} for which $G$ is 
infinite cyclic. 

These results have a bearing on the following well-known problem. 

\begin{conj}[Universal Covering Conjecture] Let $X$ be a closed,  
connected, \irr, orientable \tm\ with $\pi_1(X)$ infinite. Then the universal 
covering space $\widetilde{X}$ of $X$ is \hm\ to \RRR. \end{conj} 

Since there are only countably many homeomorphism types of closed 
3-manifolds Theorem 2 implies that there must exist uncountably many  
\rirr\ Whitehead manifolds \Mt\ which cover open 3-manifolds $M$ with 
$\pi_1(M)\cong G$ but cannot cover a closed 3-manifold. 
This generalizes a result of Tinsley and Wright \cite{Ti-Wr} which shows that 
there must exist uncountably many non-\rirr\ Whitehead manifolds 
\Mt\ which cover open 3-manifolds $M$ with $\pi_1(M)$ infinite cyclic 
but cannot cover a closed \tm. Unfortunately 
this argument does not provide any \textit{specific} such examples. 
Specific examples of non-\rirr\ Whitehead manifolds \Mt\ which cover 
open 3-manifolds $M$ with $\pi_1(M)$ infinite cyclic or, more generally, 
a countable free group, but cannot cover a closed 3-manifold are given 
in \cite{My cover} and \cite{My free}, respectively. At the time of this 
writing the problem of providing specific examples 
of \rirr\ Whitehead manifolds which non-trivially cover other open 
3-manifolds but cannot cover a closed 3-manifold is still open. 

One can make several conjectures related to Conjecture 1. We consider 
the selection below. In all of them $G$ is assumed to be a finitely 
generated group of covering translations acting on a Whitehead manifold 
$W$ with quotient a 3-manifold $M$. 

\begin{conj} $G$ is a free product of infinite cyclic groups and 
fundamental groups of \birr\ Haken manifolds. \end{conj} 

\begin{conj} $G$ is a free group or contains an infinite closed surface 
group. \end{conj}

\begin{conj} If $W$ is \rirr, then $G$ is a free product of infinite 
cyclic groups and infinite closed surface groups. \end{conj}

A proper plane $\Pi$ in $W$ is \textbf{equivariant} if for each 
$g \in G$ either $g(\Pi)=\Pi$ or $\Pi \cap g(\Pi)=\nul$. 

\begin{conj}[Special Equivariant Plane Conjecture] If $G$ is not a 
free product of infinite cyclic groups and 
infinite closed surface groups, then $W$ contains a non-trivial 
equivariant plane. \end{conj}

\begin{conj}[Equivariant Plane Conjecture] If $W$ contains a non-trivial 
plane, then it contains a non-trivial equivariant plane. \end{conj} 

These conjectures are related as follows. 

\begin{thm} $(4)\Leftarrow(1)\Leftrightarrow(2)\Leftrightarrow(3)
\Leftrightarrow(5)\Leftarrow(4+6)$ \end{thm}

Theorems 1 and 3 are proven in section 2. Theorem 2 is proven in sections 
3--7. Section 3 presents a modified version of the criterion used by Scott 
and Tucker \cite{Sc-Tu} for showing that a \tm\ is \rirr. Sections 4 and 
5 treat, respectively, the special cases in which $G$ is an infinite cyclic 
group and an infinite closed surface group. The constructions and notation of 
these special cases are used in section 6, which treats the general 
case. Section 7 shows how to get uncountably many $M$ with non-\hm\ \Mt\  
for each group $G$.

\section{The proofs of Theorems 1 and 3}

\begin{lem} Let $M$ be a connected, \pirr, open \tm. Let $Q$ be a 
compact, connected, 3-dimensional submanifold of $M$ such that $\bd Q$ 
is \inc\ in $M$ and $\pi_1(Q)$ is not an infinite closed 
surface group. Let $p:\Mt \rightarrow M$ be the universal covering map and 
$G$ the group of covering translations. Let \Qt\ be a component of $p\n(Q)$. 
Then \begin{enumerate} 
\item Each component of $p\n(\bd Q)$ is a plane. 
\item There is no component $\Pi$ of $\bd \Qt$ which is invariant under 
the subgroup $G_0$ of $G$ consisting of those covering translations which 
leave $\Qt$ invariant. 
\item If each component of $\bd \Qt$ is a trivial plane, then \Mt\ is 
\hm\ to \RRR. \end{enumerate} \end{lem}

\begin{proof} (1) follows from the incompressibility of $\bd Q$ in $M$. 

Suppose $S$ is a component of $\bd Q$ and $\Pi$ is a component of 
$p\n(S)$ which is invariant under $G_0$. Since the restriction of $p$ to 
\Qt\ is the universal covering space of $Q$ and the restriction of $G_0$ to 
\Qt\ is the group of covering translations we have that 
$\pi_1(S) \rightarrow \pi_1(Q)$ is an isomorphism, contradicting our 
assumption on $\pi_1(Q)$. This establishes (2). 

We now prove (3). Suppose that each component $\Pi$ of $\bd \Qt$ bounds an 
end-proper halfspace $H_{\Pi}$ in \Mt. Let $K_{\Pi}$ be the closure of the 
component of $\Mt-\Pi$ which does not contain $\inte \Qt$. 

Assume that for all such $\Pi$ we have $H_{\Pi}=K_{\Pi}$. Then \Mt\ is 
the union of \Qt\ and an open collar attached to $\bd \Qt$, hence \Mt\ is 
\hm\ to $\inte \Qt$. Since $Q$ is Haken, the Waldhausen Compactification 
Theorem \cite{Wa} implies that  
\Qt\ is \hm\ to a closed 3-ball minus a closed subset of its boundary, hence 
$\inte \Qt$ is \hm\ to \RRR, and we are done. 

Thus we may assume that for some $\Pi$ we have $H_{\Pi}\neq K_{\Pi}$. 
Then $H_{\Pi}\cap K_{\Pi}=\Pi$ and $H_{\Pi}\cup K_{\Pi}=\Mt$. 
Now $G_0$ has an element $g$ such that $g(\Pi)\neq \Pi$. Since 
$\Qt \subseteq H_{\Pi}$ and $g(\Qt)=\Qt$ we must have $g(K_{\Pi})
\subseteq H_{\Pi}$. Since \halfspace\ is \rirr\ (see e.g. \cite{My endsum}) 
it follows 
that $K_{\Pi}$ is \hm\ to \halfspace. Thus \Mt\ is \hm\ to \RRR. \end{proof}

\begin{proof}[Proof of Theorem 1] By passing to a covering space of 
$M$, if necessary, we may assume that $\pi_1(M)$ is indecomposable with 
respect to free products and is neither an infinite cyclic group nor an 
infinite closed surface group. Let $C$ be the Scott compact core \cite{Sc} of 
$M$, i.e. $C$ is a compact, connected, 3-dimensional submanifold of $M$ 
such that $\pi_1(C) \rightarrow \pi_1(M)$ is an isomorphism. The conditions 
on $\pi_1(M)$ imply that $\bd C$ is \inc\ in $M$. We thus can apply 
Lemma 2.1 with $Q=C$ to finish the proof. \end{proof}

\begin{proof}[Proof of Theorem 3] We first show that 
$(1)\Rightarrow(2)\Rightarrow(3)\Rightarrow(1)$. If (1) is true, then 
$M$ must be non-compact; this follows from the fact that if $M$ were closed 
and non-orientable, then it would be Haken and so have universal covering 
space \hm\ to \RRR. Let $C$ be the Scott compact core for $M$. 
Since $M$ is \irr\ we may assume that no component of $\bd C$ is a 
2-sphere; it follows that $C$ is \irr. If $C$ is \birr, then we are done. 
If $C$ is not \birr, then there is a finite set of compressing disks for 
$\bd C$ in $C$ which express $C$ as a \bd-connected sum of 3-balls and 
\birr\ Haken manifolds, thus yielding (2). Clearly $(2)\Rightarrow(3)$. 
Suppose (3) is true and $M$ is closed. If $G$ is free, then $M$ is by 
\cite[Theorem 5.2]{He} a connected sum of 2-sphere bundles over $S^1$, 
hence is not aspherical, 
hence $W$ is not contractible. If $G$ contains an infinite  closed surface 
group, then 
by a result of Hass, Rubinstein, and Scott \cite{HRS} $W$ is \hm\ to \RRR. 

Clearly Theorem 1 and the fact that $M$ cannot be closed and non-orientable 
show that $(1)\Rightarrow(4)$. 

We now show that 
$(1)\Rightarrow(5)$. Let $C$ be the Scott compact core of $M$. Then the 
assumptions on $G$ imply that there is a set of compressing disks for 
$\bd C$ in $C$ such that some component $Q$ of $C$ split along this 
collection of disks satisfies the hypotheses of Lemma 2.1. Thus any 
component of the pre-image of $\bd Q$ is an equivariant non-trivial 
plane. 

We next show that $(5)\Rightarrow(1)$. Assume $M$ is closed. If $\pi_1(M)$ 
is a free product of infinite cyclic groups and infinite closed surface 
groups, then we apply (3) to obtain (1). If $\pi_1(M)$ is not such a 
group, then the existence of an equivariant plane, together with the 
compactness of $M$, implies that $M$ is Haken, and so (1) follows by 
Waldhausen \cite{Wa}. 

Finally we show that $(4+6)\Rightarrow(1)$. If $W$ is \rirr, then (4) 
implies the hypothesis of (2), hence implies (1). If $W$ is not 
\rirr, then (6) implies as before that $M$ is Haken, thus (1) holds. 
\end{proof}

\section{Nice quasi-exhaustions and \RR-irreducibility}

We shall reformulate a criterion due to Scott and Tucker \cite{Sc-Tu} 
for a \pirr\ open \tm\ to be \rirr. A proper plane $\Pi$ in an open  \tm\ 
$W$ is \textbf{homotopically trivial} if for any compact subset $C$ of $W$ the 
inclusion map of $\Pi$ is end-properly homotopic to a map whose image 
is disjoint from $C$. 

\begin{lem} Let $W$ be an \irr, open \tm, and let $\Pi$ be a proper plane in 
$W$. If $\Pi$ is homotopically trivial, then $\Pi$ is trivial. \end{lem}

\begin{proof} This is Lemma 4.1 of \cite{Sc-Tu}. \end{proof}

\begin{lem} Let $W$ be a connected, \irr, open \tm, and let 
$\{C_n\}_{n\geq 1}$, be a 
sequence of compact 3-dimensional submanifolds of $W$ such that $C_n 
\subseteq \inte C_{n+1}$ and 
\begin{enumerate}
\item each $C_n$ is \irr, 
\item each $\bd C_n$ is \inc\ in $W-\inte C_n$, 
\item if $D$ is a proper disk in $C_{n+1}$ which is in general position 
with respect to $\bd C_n$ such that $\bd D$ is not null-homotopic in 
$\bd C_{n+1}$, then $D\cap \bd C_n$ has at least two components which 
are not null-homotopic in $\bd C_n$ and bound disjoint disks in $D$. 
\end{enumerate}
Then any proper plane in $W$ can be end-properly homotoped off 
$C_n$ for any $n$. \end{lem}

\begin{proof} This is Lemma 4.2 of \cite{Sc-Tu}. \end{proof}

The precise criterion we shall use is as follows. 

\begin{lem} Let $W$ be a connected, \irr, open \tm. Suppose that for 
each compact 
subset $K$ of $W$ there is a sequence $\{C_n\}_{n\geq 1}$ of compact 
3-dimensional submanifolds such that $C_n \subseteq \inte C_{n+1}$ and 
\begin{enumerate}
\item each $C_n$ is \irr, 
\item each $\bd C_n$ is \inc\ in $W-\inte C_n$ and has positive genus, 
\item each $C_{n+1}-\inte C_n$ is \irr, \birr, and anannular, 
\item $K \subseteq C_1$.  
\end{enumerate}
Then $W$ is \rirr. \end{lem}

\begin{proof} Let $D$ be a disk as in part (iii) of Lemma 3.2. If every 
component of $D \cap \bd C_n$ is null-homotopic in $\bd C_n$, then one can 
isotop $D$ so that $D\cap C_n=\nul$ and hence $\bd C_{n+1}$ is 
compressible in $C_{n+1}-\inte C_n$. If only one component $\alpha$ of 
$D \cap \bd C_n$ is not null-homotopic in $\bd C_n$, then $\bd D \cup 
\alpha$ bounds an annulus $A$ which can be isotoped so that $A \cap 
\bd C_n=\alpha$, hence $C_{n+1}-\inte C_n$ is not anannular. If no two 
of the components of $D \cap \bd C_n$ which are not null-homotopic in 
$\bd C_n$ bound disjoint disks in $D$, then these components must be 
nested on $D$. We can isotop $D$ to remove null-homotopic components and then 
intermediate annuli to again get an incompressible annulus joining 
$\bd C_{n+1}$ to $\bd C_n$. Now apply Lemma 3.2 and then Lemma 3.1. 
\end{proof}

Let $\{C_n\}$ be a sequence of compact, connected 3-dimensional submanifolds 
of an \irr, open \tm\ $W$ such that $C_n \subseteq \inte C_{n+1}$ such that 
$W-\inte C_n$ has no compact components. This will be called a  
\textbf{quasi-exhaustion} for $W$. A quasi-exhaustion for $W$ whose union 
is $W$ is an \textbf{exhaustion} for $W$. A quasi-exhaustion is \textbf{nice} 
if it satisfies conditions (1)--(3) of Lemma 3.3. Thus that lemma can be 
rephrased by saying that if every compact subset of $W$ is contained in 
the first term of a nice quasi-exhaustion, then $W$ is \rirr. 

We shall need some tools for constructing Whitehead manifolds with 
nice quasi-exhaustions. Define a compact, connected \tm\ $Y$ to be 
\textbf{nice} it is is \pirr, \birr, and anannular, it contains a 
two-sided proper \inc\ surface, and it is not a 3-ball; define it to 
be \textbf{excellent} if, in addition, every connected, proper, 
incompressible surface of zero Euler characteristic in $Y$ is \bd-parallel. 
So in particular an excellent \tm\ is anannular and atoroidal while a nice 
\tm\ is anannular but may contain a non-\bd-parallel \inc\ torus. 

A proper 1-manifold in a compact \tm\ is \textbf{excellent} if its exterior 
is excellent; it is \textbf{\pe\ } if the union of each non-empty subset of 
the set of its components is excellent. 

\begin{lem} Every proper 1-manifold in a compact, connected \tm\ whose 
boundary contains no 2-spheres or projective planes is homotopic rel \bd\ 
to an excellent proper 1-manifold. \end{lem}

\begin{proof} This is a special case of Theorem 1.1 of \cite{My excel}. 
\end{proof}

Define a \textbf{$k$-tangle} to be a disjoint union of $k$ proper arcs in a 
3-ball. 

\begin{lem} For all $k\geq 1$ \pe\ $k$-tangles exist. \end{lem} 

\begin{proof} This is Theorem 6.3 of \cite{My attach}. \end{proof}

We shall also need the following criterion for gluing together excellent 
3-manifolds to get an excellent 3-manifold. 

\begin{lem} Let $Y$ be a compact, connected \tm. Let $S$ be a compact, 
proper, two-sided surface in $Y$. Let $Y\p$ be the \tm\ obtained by 
splitting $Y$ along $S$. Let $S\p$ and $S\pp$ be the two copies of $S$ 
which are identified to obtain $Y$. If each component of $Y\p$ is 
excellent, $S\p$, $S\pp$, and $(\bd Y\p)-\inte (S\p \cup S\pp)$ are 
\inc\ in $Y\p$, and each component of $S$ has negative Euler characteristic, 
then $Y$ is excellent. \end{lem}

\begin{proof} This is Lemma 2.1 of \cite{My excel}. \end{proof}

\section{The infinite cyclic case}

In \cite{Sc-Tu} Scott and Tucker described an \rirr\ Whitehead manifold 
which is an infinite cyclic covering space. This section gives a general 
procedure for constructing such examples. The construction introduced here 
will be incorporated into that for the general case in section 6. 

Let $\Pn=\Dn \times [0,1]$, where \Dn\ is the disk of radius $n$. 
We call \Pn\ a \textbf{pillbox}. Identify $\Dn \times \{0\}$ with 
$\Dn \times \{1\}$ to obtain a solid torus \Qn. Let \Rn\ be a solid 
torus and \Hn\ a 1-handle $D \times [0,1]$ joining $\bd \Dn \times (0,1)$ 
to $\bd \Rn$. Let $\Vn=\Pn \cup \Hn \cup \Rn$ and $\Mn=\Qn \cup \Hn \cup 
\Rn$. We call \Vn\ an \textbf{eyebolt}. We embed \Mn\ in the interior of 
\Mnn\ as follows. 

We choose a collection of arcs \thz, \tho, \alz, \alo, \alt, \beo, \bet, 
\gao, \gat, \deo, \detw, \ep\ in \Pnn which satisfy certain conditions 
described below. \thz, \tho, and \alz\ meet in a common endpoint in 
$\inte \Pnn$ but are otherwise disjoint. The other endpoints of \thz\ and 
\alz\ lie in $(\inte \Dnn) \times \{0\}$; that of \tho\ lies in $(\inte 
\Dnn) \times \{1\}$. We let $\tht=\thz \cup \tho$. All the other arcs are 
proper arcs in \Pnn which are disjoint from each other and from 
$\tht \cup \alz$. \gao, \bet, and \detw\ run from $(\inte \Dnn)\times \{0\}$ 
to itself. \gat, \beo, and \deo run from $(\inte \Dnn)\times \{1\}$ to 
itself. \alo\ runs from $(\inte \Dnn) \times \{0\}$ to 
$(\inte \Dnn) \times \{1\}$. \alt\ runs 
from $(\inte \Dnn) \times \{1\}$ to $\inte (\Pnn \cap \Hnn)$. \ep\ runs from 
$\inte (\Pnn \cap \Hnn)$ to itself. We denote the image in \Qnn\ of an 
arc by the same symbol, relying on the context to distinguish an arc in 
\Pnn\ from its image in \Qnn. We require that \tht\ be a simple closed curve 
in \Qnn\ and that $\alz \cup \beo \cup \gao \cup \deo \cup \alo \cup 
\bet \cup \gat \cup \detw$ is an arc consisting of subarcs which occur in the 
given order. We require that any union of these arcs which contains \alz\ 
and at least one other arc has excellent exterior in \Pnn, and that the 
same is true for any union of these arcs which contains neither \thz, \tho, 
nor \alz. This can be achieved as follows. Note that the exterior of \alz\ 
in \Pnn\ is a 3-ball $B$. Choose a poly-excellent 11-tangle in $B$ and then 
slide its endpoints so that exactly two of the arcs meet a regular 
neighborhood of \alz. Extend them to meet \alz\ in the desired configuration. 

Next let \kao, \kat, and \kad\ be product arcs in \Hnn\ joining 
$(\inte D) \times \{0\}$ to $(\inte D) \times \{1\}$. Let 
$\Rn \subseteq \inte \Rnn$ be any null-homotopic embedding. Let \lao\ 
and \lat\ be disjoint proper arcs in $\Rnn-\inte \Rn$ with \lao\ 
joining $\inte (\Hnn \cap \Rnn)$ to itself and \lat\ joining 
$\inte (\Hnn \cap \Rnn)$ to $\bd \Rn$. We require $\lao \cup \lat$ 
to be excellent in $\Rnn-\inte \Rn$. We also require that these arcs, 
together with \ep, fit into an arc whose subarcs form the sequence 
\kao, \lao, \kat, \ep, \kad, \lat\ and that \kao\ meets \alt\ in a 
common endpoint. 

Now we embed \Pn\ in \Pnn\ as a regular neighborhood of the arc \tht\ so that 
the two disks of $\Pn \cap (\Dnn \times \{0,1\})$ are identified to give 
an embedding of \Qn\ in \Qnn. Note that these embeddings are not consistent 
with the product structures. From the discussion above we have an arc \om\ 
in $\Mnn-\inte(\Qn \cup \Rn)$ running from $\bd \Qn$ to $\bd \Rn$. We 
embed \Hn\ as a regular neighborhood of \om. We change notation slightly 
by now letting \alz\ be the old \alz\ minus its intersection with the 
interior of \Qn. 

We let $M$ be the direct limit of the \Mn\ and let $p:\Mt \rightarrow M$ 
be the universal covering map. $p\n(\Qn)=p\n(\Pn)$ is the union of pillboxes 
$\Pnj=\Dn \times [j,j+1]$ meeting along the $\Dn \times \{j\}$ to form 
$\Dn \times \R$. Note that this embedding is not the product embedding. 
$p\n(\Rn)$ is a disjoint union of solid tori \Rnj. $p\n(\Hn)$ is a disjoint 
union of 1-handles \Hnj\ joining $\bd \Dn \times (j,j+1)$ to $\bd \Rnj$; 
these are regular neighborhoods of lifts \omj\ of \om. $p\n(\Mn)=
p\n(\Vn)$ is the union of $p\n(\Pn)$, $p\n(\Hn)$, and $p\n(\Rn)$. It is 
the union of eyebolts $\Vnj=\Pnj \cup \Hnj \cup \Rnj$ meeting along the 
$\Dn \times \{j\}$. \Mt\ is the nested union of the $p\n(\Mn)$. 

Let $\Simn=\cup^m_{j=-m} \Vnj$ and $\Lamn=P_{n,-(m+1)} \cup P_{n,m+1}$. 
Let $\Phi^m_1=\nul$, and, for $n\geq 2$,
let $\Phmn=\cup^{m+n}_{j=m+2} (P_{n,-j} \cup \Pnj)$.   
Note that \Lamn\ and \Phmn\ (for $n\geq 2$) are each 
disjoint unions of two 3-balls, $\Lamn \cap \Simn$ is a pair of disjoint 
disks, and (for $n\geq 2$) so is $\Lamn \cap \Phmn$. Define $\Cmn=
\Simn \cup \Lamn \cup \Phmn$. 

\begin{lem} $\{C^m_m\}$ is an exhaustion for \Mt. Each $C^m$ is a nice 
quasi-exhaustion. \end{lem} 

\begin{proof} $\Cmn \subseteq \inte \Cmnn$, and $\Cmn \subseteq C^{m+1}_n$. 
A given compact subset $K$ of \Mt\ lies in some $p\n(\Mn)$ and thus in a 
finite union of \Vnj\ and hence in some $\Simn \subseteq \Cmn \subseteq 
C^q_q$, where $q=\max\{m,n\}$. Thus $\{C^m_m\}$ is an exhaustion for \Mt. 

\Cmn\ is a cube with $2m+1$ handles. Let $Y=\Cmnn-\inte \Cmn$. We will show 
that $Y$ is excellent by successive applications of Lemma 3.6. 

Consider a \Pnnj\ contained in \Cmnn. If $|j|<m$, then it meets \Cmn\ in 
a regular neighborhood of the union of the $j^{th}$ copies of all the arcs in 
\Pnn. Thus $Y \cap \Pnnj$ is excellent, and Lemma 3.6 implies that the union 
of these $Y \cap \Pnnj$ is excellent. For $|j|\geq m$ some care must be 
taken so that one is always gluing excellent 3-manifolds along surfaces of 
the appropriate type. Note that $Y \cap (P_{n+1,m} \cup P_{n+1,m+1} 
\cup \cdots \cup P_{n+1,m+n-1} \cup P_{n+1,m+n})$ is equal to the exterior 
of the $m^{th}$ copy of all the arcs but \beo\ and \detw\ in $P_{n+1,m}$ 
together with the exterior of the $(m+1)^{st}$ copy of \bet, \detw, and \tht\ 
in $P_{n+1,m+1}$, the exterior of the $j^{th}$ copy of \tht\ in \Pnnj\ for 
$m+1 < j < m+n$, and the 3-ball $P_{n+1,m+n}$. This space is \hm\ to the 
exterior of the $m^{th}$ copy of all the arcs but \beo, \deo, and \tho\ in 
$P_{n+1,m+1}$ together with the exterior of the $(m+1)^{st}$ copy of \bet\ 
and \detw\ in $P_{n+1,m+1}$, and the 3-ball consisting of the union of the 
\Pnnj\ for which $m+1<j\leq m+n$. This can be seen by taking the arc 
consisting of the $m^{th}$ copy of \tho\ and the $j^{th}$ copy of \tht\ 
for $m<j<m+n$ and retracting it onto the endpoint in which it meets the 
rest of the graph. This space is then excellent by Lemma 3.6. 
Similar remarks apply for $j\leq -m$, so these spaces can be added on to 
get that $Y \cap \cup^{m+n}_{j=-(m+n)} \Pnnj$ is excellent. 

We fill in the remainder of $Y$ by adding the exteriors of the $j^{th}$ 
copies of \kao, \kat, and \kad\ in \Hnnj\ and $\lao \cup \lat$ in 
$\Rnnj-\inte \Rnj$ for $|j|\leq m$. Since the first of these spaces is a 
product the union of the two spaces is \hm\ to the second space, and  
Lemma 3.6 applies to complete the proof that $Y$ is excellent.

It remains to show that each $\bd \Cmn$ is \inc\ in $\Mt-\inte \Cmn$. 
Since each $C^m_{n+s+1}-\inte C^m_{m+s}$ is \birr\ we have that 
$\bd \Cmn$ is \inc\ in $C^m_{n+q}-\inte \Cmn$ for each $q\geq 1$. 
$p\n(M_{n+q})$ is the union of $C^m_{n+q}$ and the sets $C^-=\cup_{j<-m} 
V_{n+q,j}$ and $C^+=\cup_{j>m} V_{n+q,j}$. We have that $C^- \cap C^m_{n+q}$ 
and $C^+ \cap C^m_{n+q}$ are disjoint disks, while $C^- \cap C^+=\nul$. 
It follows that $\bd \Cmn$ is \inc\ in $p\n(M_{n+q})-\inte \Cmn$. 
Since \Mt\ is the nested union of the $p\n(M_{n+q})$ over all $q\geq 1$ 
we have the desired result. \end{proof}

\section{The surface group case}

Let $F$ be a closed, connected surface other than $S^2$ or \PP. Let 
$n\geq 1$. Regard $F$ as being obtained from a $2k$-gon $E$, $k\geq 2$, 
by identifying sides $s_i$ and $s_i\p$, $1\leq i \leq k$. This induces 
an identification of the lateral sides $S_i=s_i \times [-n,n]$ and 
$S_i\p=s_i \times [-n,n]$ of the \textbf{prism} $\Pn=E \times [-n,n]$ 
which yields $\Qn=F \times [-n,n]$. Let $\Rn$ be a solid torus and $\Hn$ 
a 1-handle $D \times [0,1]$. Let $\Vn=\Pn \cup \Hn \cup \Rn$, where 
$\Pn \cap \Rn=D \times \{1\}$ is a disk in $\bd \Rn$, and 
$\Pn \cap \Rn=\nul$. We again call \Vn\ an \textbf{eyebolt}. It is a 
solid torus whose image under the identification is $\Mn=\Qn \cup \Hn 
\cup \Rn$, a space \hm\ to the \bd-connected sum of $F\times [-n,n]$ and 
a solid torus. 

We define an open \tm\ $M$ by specifying an embedding of \Mn\ in the 
interior of \Mnn\ and letting $M$ be the direct limit. The inclusion 
$[-n,n] \subseteq [-(n+1),n+1]$ induces $\Pn \subseteq \Pnn$ and hence 
$\Qn \subseteq \Qnn$. We let $\Rn \subseteq \inte \Rnn$ be any 
null-homotopic embedding. Again the interesting part of the embedding 
will be that of \Hn\ in \Mnn. It will be the regular neighborhood of a 
certain arc \om\ in $\Mnn-\inte(\Qn \cup \Rn)$ joining $\bd \Qn$ to 
$\bd \Rn$. 

The arc \om\ is the union of $4k+7$ arcs any two of which are either 
disjoint or have one common endpoint. The $4k+2$ arcs \alz, \ali, \bei, 
\gai, \dei, $1\leq i \leq k$, and \ep\ lie in $E \times [n,n+1]$ and 
are identified with their images in \Qnn; the three arcs \kao, \kat, 
and \kad\ lie in \Hnn, and the two arcs \lao\ and \lat\ lie in \Rnn. 
These arcs will have special properties to be described later. We first 
describe their combinatorics.  The arcs in \Pnn\ are all proper arcs in 
$E \times [n,n+1]$. \alz\ runs from $(\inte E) \times \{n\}$ to $\inte S_1$. 
For $1 \leq i < k$, \ali\ runs from $\inte S_i$ to $\inte S_{i+1}$. 
\alk\ runs from $\inte S_k\p$ to $\inte (\Pnn \cap \Hnn)$. 
For $1 \leq i \leq k$, \bei\ and \dei\ each run from $\inte S_i\p$ to 
itself, while \gai\ runs from $\inte S_i$ to itself. These arcs are chosen 
so that under the identification their endpoints match up in such a way 
as to give  an arc which follows the sequence \alz, \beo, \gao, \deo, \alo, 
\ldots, \bek, \gak, \dek, \alk. We require \ep\ to run from 
$\inte(\Pnn \cap \Hnn)$ to itself. \kao, \kat, and \kad\ are product arcs 
in \Hnn\ lying in $(\inte D) \times [0,1]$. \lao\ and \lat\ are proper arcs 
in $\Rnn-\inte \Rn$, with \lao\ running from $\inte (\Hnn \cap \Rnn)$ to 
itself and \lat\ running from $\inte(\Hnn \cap \Rnn)$ to $\bd \Rn$. These 
arcs are chosen so as to fit together into the sequence \kao, \lao, \kat, 
\ep, \kad, \lat\ with the endpoint of \kao\ other than $\kao \cap \lao$ 
being the same as the endpoint of \alk\ other than $\alk \cap \dek$. This 
gives \om. 

We now describe the special properties required of these arcs. We require 
that $\alz \cup \beo \cup \gao \cup \deo \cup \alo \cup \cdots \cup \bek \cup 
\gak \cup \dek \cup \ep$ be a \pe\ $(4k+2)$-tangle in $E \times [n,n+1]$ 
and $\lao \cup \lat$ to be an excellent 1-manifold in $\Rnn-\inte \Rn$. 

We now consider the universal covering map $p:\Mt \rightarrow M$. Our goal 
is to construct a sequence $\{C^m\}$ of nice quasi-exhaustions whose 
diagonal $\{C^m_m\}$ is an exhaustion for \Mt. 

The universal covering space \Ft\ of $F$ is tesselated by copies \Ej\ of $E$. 
We fix one such copy $E_1$. We inductively define an exhaustion $\{\Fm\}$ 
for \Ft\ as follows. $F_1=E_1$. \Fmm\ is the union of \Fm\ and all those 
\Ej\ which meet it. Each \Fm\ is a disk (which we call a \textbf{star}). 
The \textbf{inner corona} \Icm\ of \Fm\ is the annulus $\Fmm-\inte \Fm$. 
Each vertex on $\bd \Fm$ lies in either one or two of those \Ej\ contained 
in \Fm. Each \Ej\ in \Icm\ meets meets \Fm\ in either an edge or a vertex; 
in both cases it meets exactly two adjacent $E_\ell$ of $\Icm$, and each of 
these intersections is an edge. For $n\geq 2$ we define the \textbf{outer 
$n$-corona} \Omn\ to be the annulus $F_{m+n}-\inte \Fmm$; we define 
$O^m_1=\nul$. 
Let $\sigma_2$ be a proper arc in $F_2$ consisting of three edges of the 
polygons in $F_2$. . Inductively define a proper arc $\sigma_{m+1}$ in 
\Fmm\ by adjoining to $\sigma_m$ two arcs spanning $\Icm$ which are edges 
of polygons in \Icm. Thus each $\sigma_m$ is an edge path in \Fm\ splitting 
it into two unions of polygons \Fmp\ and \Fmpp. 

We now consider the structure of \Mt. For $n\geq 1$, $p\n(\Qn)=p\n(\Pn)$ is 
the union of prisms $\Pnj=\Ej \times [-n,n]$ meeting along their lateral 
sides to form $\Ft \times [-n,n]$. $p\n(\Rn)$ is a disjoint union of solid 
tori \Rnj. $p\n(\Hn)$ is a disjoint union of 1-handles \Hnj\ running from 
$\Ej\times \{n\}$ to $\bd \Rnj$; these are regular neighborhoods of lifts 
\omj\ of \om. Now $p\n(\Mn)=p\n(\Vn)$ is the union of $p\n(\Pn)$, $p\n(\Hn)$, 
and $p\n(\Rn)$. It can be expressed as the union of the eyebolts $\Vnj=
\Pnj \cup \Hnj \cup \Rnj$ meeting along the lateral sides of the \Pnj. 
Finally \Mt\ is the nested union of the $p\n(\Mn)$. 

Let \Simn\ be the union of those \Vnj\ such that \Ej\ is in the star \Fm. 
Let \Lamn\ be the union of those \Pnj\ such that \Ej\ is in the inner corona 
\Icm. Let \Phmn\ be the union of those \Pnj\ such that \Ej\ is in the 
outer $n$-corona \Omn. Note that \Lamn\ and \Phmn\ (for $n\geq 2$) are solid
tori, $\Lamn \cap \Simn$ is an annulus which goes around $\Lamn$ once 
longitudinally and consists of those lateral sides of the prisms in \Simn\ 
which lie on $\bd \Simn$, and (for $n \geq 2$) $\Lamn \cap \Phmn$ is an 
annulus which goes around each of these solid tori once longitudinally. 
We now define $\Cmn=\Simn \cup \Lamn \cup \Phmn$. 

\begin{lem} $\{C^m_m\}$ is an exhaustion for \Mt. Each $C^m$ is a nice 
quasi-exhaustion. \end{lem}

\begin{proof} Note that $\Cmn \subseteq \inte \Cmnn$, and $\Cmn \subseteq 
\Cmmn$. Suppose $K$ is some compact subset of \Mt. Then $K$ lies in 
some $p\n(\Mn)$ and thus in a finite union of \Vnj\ and hence in some 
$\Simn \subseteq \Cmn \subseteq C^q_q$, where $q=\max\{m,n\}$. Thus 
$\{C^m_m\}$ is an exhaustion for \Mt. 

Each \Cmn\ is a cube with handles, so is \irr. The number of handles is at 
least one, so $\bd \Cmn$ has positive genus. Let $Y=\Cmnn-\inte \Cmn$. 
We will prove that $Y$ is excellent by successive applications of Lemma 3.6. 
Let $\Pnnjp$ and $P^-_{n+1,j}$ denote, respectively, 
$\Ej \times [n,n+1]$ and $\Ej \times [-(n+1),-n]$.

Consider a \Pnnj\ contained in \Simnn. It meets \Cmn\ in \Pnnj\ 
together with regular neighborhoods of certain arcs in \Pnnjp. These 
arcs consist at least of the $j^{th}$ copies of the \ali, the \gai, and 
\ep\ which are part of the lift \omj\ of \om. If another prism $P_{n+1,\ell}$ 
in \Simnn\ meets \Pnnj\ in a common lateral side, then either \omj\ 
or $\om_{\ell}$ will meet this side; in the latter case this contributes a 
\bei\ and \dei\ to the subsystem of arcs in \Pnnjp. Since the full system 
of arcs was chosen to be \pe\ this subsystem of arcs is excellent and so has 
excellent exterior $Y \cap \Pnnjp$. Let $U\p$ be the union of those 
$Y\cap \Pnnjp$ such that $\Ej \subseteq \Fmp$. This space can be built up 
inductively by gluing on one $Y \cap \Pnnjp$ at a time, with the gluing 
being done along either a disk with two holes (when \Pnnj\ is glued along 
one lateral side) or a disk with four holes (when \Pnnj\ is glued along 
two adjacent lateral sides). No component of the complement of this surface 
in the boundary of either manifold is a disk, hence this surface is \inc\ 
in each manifold. It follows that $U\p$ is excellent. Similar remarks apply 
to the space $U\pp$ associated with $\Fmpp$. 

Next consider a \Pnnjp\ contained in \Lamnn. If $\Ej \subseteq \Fmm$ and 
meets \Fmp\ in an edge of $E_{\ell} \subseteq \Fmp$, then either $\om_{\ell}$ 
misses \Pnnjp\ or meets it in copies of \bei\ and \dei. Thus enlarging 
$U\p$ by adding $Y \cap \Pnnjp$ either adds a 3-ball along a disk in its 
boundary, giving a space \hm\ to $U\p$ or gives a new excellent \tm. 
We adjoin all such $Y\cap \Pnnjp$ to $U\p$. Then we consider those $\Ej$ 
which meet \Fmp\ in a vertex. Then $\Pnnjp=Y\cap\Pnnjp$, and one can 
successively adjoin these 3-balls along disks in their boundaries. We 
denote the enlargement of $U\p$ from all these additions again by $U\p$. 
Similar remarks apply to $U\pp$. 

Now $(F\p_{m+n+1}-\inte \Fmmp) \times [n,n+1]$ is a 3-ball which meets 
$U\p$ in a disk, so we adjoin it to $U\p$ to get a new $U\p$ \hm\ to the 
old one. We then adjoin the 3-ball $(F\p_{m+n+1}-\inte F\p_{m+n}) \times 
[-n,n]$ which meets this space along a disk to obtain our final $U\p$. 
The same construction gives $U\pp$. 

Now $U\p$ and $U\pp$ are each excellent. $U\p \cap U\pp$ is an annulus with a 
positive number of disks removed from its interior corresponding to its 
intersection with arcs passing from $\Fmp \times [n,n+1]$ to 
$\Fmpp \times [n,n+1]$. No component of the complement of this surface in 
$\bd U\p$ or in $\bd U\pp$ is a disk; this corresponds to the fact that 
$\Fmp \times \{n\}$, $\Fmpp \times \{n\}$, $\Fmp \times \{n+1\}$, and 
$\Fmpp \times \{n+1\}$ each meet some \omj. Thus this surface is \inc\ in 
both $U\p$ and $U\pp$, so $U\p \cup U\pp$ is excellent. 

Finally we add on the $Y\cap\Rnj$ for $\Ej \subseteq \Fm$ to $U\p \cup U\pp$ 
to conclude that $Y$ is excellent. 

It remains to show that each $\bd \Cmn$ is \inc\ in $\Mt-\inte \Cmn$. 
First note that since each $C^m_{n+s+1}-\inte C^m_{n+s}$ is \birr\ 
we must have that $\bd \Cmn$ is \inc\ in $C^m_{n+q}-\inte \Cmn$ for each 
$q\geq 1$. Now consider the set 
\[ \widetilde{M}_{n+q}=p\n(M_{n+q}) \cup (\Ft \times [-(n+q+1),-(n+q)]). \] 
It can be obtained from $C^m_{n+q}$ as 
follows. First add the solid tori $R_{n+1,j} \cup H_{n+q,j}$ in 
$p\n(M_{n+q})$ for which $\Ej \subseteq F_{m+q+n}$; these meet $C^m_{n+q}$ 
in disks. Then add 
\[ (F_{m+q+n} \times [-(n+q+1),-(n+q)]) \cup 
(\Ft-\inte F_{m+q+n})\times[-(n+q+1),n+q]). \] 
This is a space \hm\ to 
\prodspace\ which meets $C^m_{n+q}$ in the disk 
\[ (F_{m+q+n} \times \{-(n+q)\})\cup((\bd F_{m+q+n})\times[-(n+q),n+q]). \] 
Lastly add all the 
remaining solid tori $R_{n+q,j} \cup H_{n+q,j}$, where $\Ej \subseteq 
\Ft-\inte F_{m+q+n}$; these do not meet $C^m_{n+q}$. This description  
shows that $C^m_{n+q} \cap (\Mt_{n+q}-\inte C^m_{n+q})$ consists of 
(finitely many) disjoint disks, and therefore $\bd \Cmn$ is \inc\ in 
$\Mt_{n+q}-\inte \Cmn$. Finally since \Mt\ is the nested union of the 
$\Mt_{n+q}$ over all $q\geq 1$ we have that $\bd \Cmn$ is \inc\ in 
$\Mt-\inte \Cmn$. \end{proof}

\section{The general case}

Suppose $G_1$, $\ldots$, $G_k$ are infinite cyclic groups and infinite 
closed surface groups. For $i=1, \ldots, k$ let \Pin\ be a pillbox or 
prism, as appropriate, with quotient \Qin\ a solid torus or product 
$I$-bundle over a closed surface, respectively. We let \Hin\ be a 
1-handle attached to \Pin\ as before. We let \Rn\ be a common solid torus 
to which we attach the other ends of all the \Hin. The union of the \Qin\ 
and \Hin\ with \Rn\ is called \Mn. As before we choose arcs in the \Pin, 
\Hin, and \Rn\ and use them to define an embedding of \Mn\ into the interior 
of \Mnn. 

The choice of arcs in $\Rnn-\inte \Rn$, as well as the embedding $\Rn 
\subseteq \inte \Rnn$, requires some discussion, since we will want this 
family $\lambda$ of arcs to be \pe. Choose a \pe\ $(2k+2)$-tangle 
$\lambda^+$ in a 3-ball $B$, with components \lait, $1\leq i\leq k+1$, 
$t=1,2$. Construct a graph in $B$ by sliding one endpoint of each 
$\lambda^i_2$, $1\leq i \leq k$, so that it lies on $\inte \lambda^{k+1}_2$. 
Thus these $\lambda^i_2$ now join $\bd B$ to distinct points on $\inte  
\lambda^{k+1}_2$; all the other \lait\ still join $\bd B$ to itself. Now 
choose disjoint disks $E_1$ and $E_2$ in $\bd B$ such that $E_t$ meets the 
graph in $\bd \lambda^{k+1}_t \cap \inte E_t$. Glue $E_1$ to $E_2$ so that $B$ 
becomes a solid torus \Rnn\ and $\lambda^{k+1}_1 \cup \lambda^{k+1}_2$ 
becomes a simple closed curve. The regular neighborhood of this simple 
closed curve is our embedding of \Rn\ in the interior of \Rnn. Clearly 
\Rn\ is null-homotopic in \Rnn. By Lemma 3.6 its exterior is excellent as is 
the exterior of the union of \Rn\ with any of the \lait, $1\leq i\leq k$, 
$t=1,2$. 

Let $p:\Mt\rightarrow M$ be the universal covering map. Then $p\n(\Rn)$ 
consists of disjoint solid tori whose union separates $p\n(\Mn)$ into 
components with closures \Limun, where \Limun\ is a component of 
$p\n(\Qin \cup \Hin)$. Let \Zimun\ be the union of \Limun\ and all those 
components of $p\n(\Rn)$ which meet it. Then $\Zimu=\cup_{n\geq 1}\Zimun$ 
is an open subset of \Mt\ which has a family $\{\Cimum\}$ of 
quasi-exhaustions as previously described. We will develop from these 
families an appropriate family $\{C^m\}$ of quasi-exhaustions of $\Mt$. 

We start by choosing a component $\widehat{R}_1$ of $p\n(R_1)$. For each 
$n$ there is then a unique component $\widehat{R}_n$ of $p\n(\Rn)$ 
which contains $\widehat{R}_1$. We define $C^1_n$ to be the union of 
$\widehat{R}_n$ and the (finitely many)  $C^{i,\mu,1}_n$ which contain it; 
this space is a solid torus which meets each of these $C^{i,\mu,1}_n$ 
in a 3-ball which is either a pillbox and a 1-handle or a prism and a 
1-handle. Suppose \Cmn\ has been defined and that it is the union of the 
\Cimumn\ for which $\Cmn \cap \Limun \neq \nul$. We define \Cmmn\ in two 
steps. We first take the union $C\p$ of all the $C^{i,\mu,m}_n$ such that 
$\Cimumn \subseteq \Cmn$. This is just the union of the $n^{th}$ elements of 
the $(m+1)^{st}$ quasi-exhaustions for those \Zimun\ such that $\{i,\mu\}$ is 
in the current index set. The second step is to enlarge the index set by 
adding those $\{i,\nu\}$ for which $C\p \cap L^{i,\nu}_n \neq \nul$ and 
then adjoin the $C^{i,\nu,m+1}_n$ to $C\p$ in order to obtain \Cmmn. One 
can observe that the \Limun\ and $p\n(\Rn)$ give $p\n(\Mn)$ a tree-like 
structure and that the passage from \Cmn\ to \Cmmn\ goes out further along 
this tree. 

\begin{lem} $\{C^m_m\}$ is an exhaustion for \Mt. $C^m$ is a nice 
quasi-exhaustion. \end{lem} 

\begin{proof} Again we have $\Cmn \subseteq \inte \Cmnn$ and $\Cmn \subseteq 
\Cmmn$ with the result that $\{C^m_m\}$ is an exhaustion for \Mt. 

As regards the excellence of $\Cmnn-\inte \Cmn$ we note that the only thing  
new takes place in those components of $p\n(\Rnn)$ contained in \Cmnn. 
Instead of two arcs $\lambda_1$ and $\lambda_2$ as before we have 
$\lambda^i_1$ and $\lambda^i_2$ as $i$ ranges over some non-empty subset of 
$\{1, \ldots, k\}$ We then apply the poly-excellence of the full set of 
\lait. 

The incompressibility of $\bd \Cmn$ in $\Mt-\inte \Cmn$ follows as 
before. We first note that $\bd \Cmn$ is \inc\ in $C^m_{n+q}-\inte \Cmn$ 
for each $q\geq 1$. Now define $\Mt_{n+q}$ to be the union of $p\n(M_{n+q})$ 
and, for each of the surface group factors $G_i$ of $G$, the copy 
$\Ft^{i,\mu} \times [-(n+q+1),-(n+q)]$ of $\Ft^i \times [-(n+q+1),-(n+q)]$ 
contained in 
\Zimu, where $\Ft^i$ is the universal covering space of the surface 
$F^i$ with $\pi_1(F^i)\cong G_i$. Then the exterior of $C^m_{n+q}$ in 
$\Mt_{n+q}$ meets it in a collection of disjoint disks, from which it 
follows that $\bd \Cmn$ is \inc\ in $\Mt_{n+q}-\inte \Cmn$, thus is \inc\ 
in $\Mt-\inte \Cmn$. \end{proof}

\section{Uncountably many examples}

We now describe how to get uncountably many examples for a given group $G$. 
We will use a trick introduced in \cite{My attach}. Let $\{X_{n,s}\}$ be a 
family of 
exteriors of non-trivial knots in $S^3$ indexed by $n\geq 2$ and $s\in
\{0,1\}$; they are chosen to be anannular, atoroidal, and pairwise 
non-\hm. (One such family is that of non-trivial, non-trefoil twist knots.) 
One chooses a function $\varphi(n)$ with values in $\{0,1\}$, i.e. a 
sequence of 0's and 1's indexed by $n$, and constructs a \tm\ $M[\varphi]$ 
by embedding $X_{n,\varphi(n)}$ in $\Mn-\inte \Mnnn$ so that 
$\bd X_{n,\varphi(n)}$ 
in \inc\ in $\Mn-\inte \Mnnn$ (but is compressible in \Mn). The idea is to 
do this in such a way that for ``large'' compact sets $C$ in \Mtphi\ one has 
components of $p\n(X_{n,\varphi(n)})$ which lie in $\Mt-C$ and have \inc\ 
boundary in $\Mt-C$ for ``large'' values of $n$; moreover, every knot 
exterior having these properties should be \hm\ to some $X_{n,\varphi(n)}$. 
Thus if \Mtphi\ and \Mtpsi\ are \hm\ one must have $\varphi(n)=\psi(n)$ 
for ``large'' $n$. One then notes that there are uncountably many functions 
which are pairwise inequivalent under this relation. 

We proceed to the details. First assume $\varphi$ is fixed, so we can write 
$s=\varphi(n)$. The most innocuous place to embed $\Xns$ is in 
$\Rn-\inte \Rnnn$ 
since this space is common to all our constructions. Recall that this space 
contains arcs $\lambda_1$, $\lambda_2$ or, if $G$ is a non-trivial free 
product, arcs $\lambda^i_1$, $\lambda^i_2$, $1\leq i\leq k$; call this 
collection of arcs $\lambda$. We wish \Xns\ to lie in the complement of 
$\lambda$ in such a way that it is \pe\ in $\Rn-\inte(\Rnnn \cup \Xns)$. 
We revise the construction of $\lambda$ from section 6 as follows. 
Let $B_0$ and $B_1$ be 3-balls. 
Choose disjoint disks $D_r$ and $D_r\p$ in $\bd B_r$. 
Let $\zeta_r$ be a simple closed curve in $\bd B_r-(D_r \cup D_r\p)$ 
which separates $D_r$ from $D_r\p$. Let $A_r$ and $A_r\p$ be the annuli 
into which $\zeta_r$ splits the annulus $\bd B_r-\inte (D_r \cup D_r\p)$, 
with the notation chosen so that $A_r \cap D_r=\nul$. Let $\tau_r$ be a 
poly-excellent $(4k+4)$-tangle in $B_r$ which is the union of 
$(2k+2)$-tangles $\rho_r$ and $\rho_r\p$ satisfying the following 
conditions. Each component of $\rho_0$ runs from $\inte D_0$ to $\inte A_0\p$. 
Each component of $\rho_1$ runs from $\inte A_1\p$ to $\inte D_1\p$. 
Each component of $\rho_0\p$ runs from $\inte D_0\p$ to itself. 
Each component of $\rho_1\p$ runs from $\inte D_1\p$ to $\inte D_1$. 
We then glue $A_0\p$ to $A_1\p$ and $D_0\p$ to $D_1\p$ so as to obtain 
a space \hm\ to a 3-ball minus the interior of an unknotted solid torus 
contained in the interior of the 3-ball. The 2-sphere boundary component 
is $D_0 \cup D_1$; the torus boundary component is $A_0 \cup A_1$. 
The gluing is done so that the endpoints of the arcs match up to give 
a system $\lambda^+$ of $2k+2$ arcs. Each arc in this system consists of 
an arc of $\rho_0$ followed by an arc of $\rho_1$ followed by an arc 
of $\rho_0\p$ followed by an arc of $\rho_1\p$. We then glue \Xns\ to 
this space along their torus boundaries so as to obtain a 3-ball $B$. 
We then apply the construction of section 6 to $\lambda^+$ to get a \pe\ 
system $\lambda$ of arcs in $\Rn-\inte \Rnnn$. It is easily seen that 
this \tm\ is nice and that $\bd \Xns$ is, up to isotopy, the unique 
\inc\ non-\bd-parallel torus in it; $\bd \Xns$ is also, up to isotopy, 
the unique \inc\ torus in the exterior \Ks\ of any non-empty union 
$\sigma$ of components of $\lambda$. 

\begin{lem} If \Mtphi\ and \Mtpsi\ are \hm\ then there is an index $N$ 
such that $\varphi(n)=\psi(n)$ for all $n\geq N$. \end{lem}

\begin{proof} Consider \Mt. 
$Y=\Cmn-\inte \Cmnnn$ contains copies of \Ks\ for various 
choices of $\sigma$. The closure of the complement in $Y$ of these copies 
consists of excellent 3-manifolds which meet the copies along \inc\ 
planar surfaces. It follows that the various copies of $\bd \Xns$ in $Y$ are, 
up to isotopy and for $n\ge3$, the unique \inc\ tori in $Y$. The 
incompressibility of $\bd \Cmn$ in $\Mt-\inte \Cmn$ implies that these 
tori are also \inc\ in $\Mt-\inte \Cmnnn$. 

Suppose $T$ is an \inc\ torus in $\Mt-\inte \Cmnnn$. Then $T$ lies in 
$\Mt_{n+q}$ for some $q \geq 0$. The exterior of $C^m_{n+q}$ in $\Mt_{n+q}$ 
consists of copies of $D \times \R$ and \prodspace\ to which disjoint 
1-handles have been attached. It meets $C^m_{n+q}$ in a set of disjoint 
disks. It follows that $T$ can be isotoped into $C^m_{n+q}-\inte \Cmnnn$. 
Since $\bd C^m_{n+u}$ for $1\leq u<q$ is not a torus it is easily seen 
that $T$ can be isotoped into some $C^m_v-\inte C^m_{v-1}$ and thus is 
isotopic to some copy of $\bd X_{v,\varphi(v)}$. Thus any knot exterior 
$X$ incompressibly embedded in $\Mt-\inte \Cmnnn$ is \hm\ 
to some $X_{v,\varphi(v)}$. 

Now consider two different functions $\varphi$ and $\psi$. We will show that 
if \Mtphi\ and \Mtpsi\ are \hm\, then there is an $N$ such that $\varphi(n)=
\psi(n)$ for all $n\geq N$. Suppose $h:\Mtphi \rightarrow \Mtpsi$ is a 
homeomorphism. Distinguish the various submanifolds arising in the 
construction of these two manifolds by appending $[\varphi]$ and $[\psi]$, 
respectively. For $n\geq 2$ there are incompressibly embedded copies 
$\widetilde{X}_{n,\varphi(n)}$ of $X_{n,\varphi(n)}$ in 
$\Mtphi-\inte C^1_1[\varphi]$. 
There is an index $\ell$ such that $h(C^1_1[\varphi]) \subseteq \inte 
C^{\ell}_{\ell}[\psi]$. By construction $\cup_{n\geq 2}
\widetilde{X}_{n,\varphi(n)}$ is end-proper in \Mtphi, so there is an index 
$N$ such that for all $n\geq N$ we have $h(\widetilde{X}_{n,\varphi(n)}) 
\subseteq \Mtpsi-\inte C^{\ell}_{\ell}[\psi]$. Since 
$h(\bd \widetilde{X}_{n,\varphi(n)})$ is \inc\ in 
$\Mtpsi-\inte h(C^1_1[\psi])$ 
it is \inc\ in the smaller set $\Mtpsi-\inte C^{\ell}_{\ell}[\psi]$. 
Thus it is \hm\ to $X_{v,\psi(v)}$ for some $v>\ell$. Since the knot 
exteriors are pairwise non-\hm\ we must have $n=v$ and $\varphi(n)=\psi(v)=
\psi(n)$. \end{proof}



\begin{thebibliography}{99}

\bibitem{HRS}
J. Hass, H. Rubinstein, and P. Scott, {\it Compactifying coverings of closed
3-manifolds}, J. Differential Geometry 30 (1989), 817--832.


\bibitem{He}
J. Hempel, {\it 3-Manifolds}, Ann. of Math. Studies, No. 86, Princeton,
(1976).

\bibitem{Ja}
W. Jaco, {\it Lectures on three-manifold topology}, CBMS Regional
Conference Series in Math., No. 43, Amer. Math. Soc. (1980).

\bibitem{Mc}
D. R. McMillan, Jr., {\it Some contractible open 3-manifolds}, Trans. 
Amer. Math. Soc. 102 (1962), 373--382. 

\bibitem{MSY}
W. Meeks, L. Simon, S. T. Yau, {\it Embedded minimal surfaces, exotic 
spheres, and manifolds with postive Ricci curvature}, Annals of Math, 
116 (1982), 621--659.

\bibitem{My excel}
R. Myers, {\it Excellent 1-manifolds in compact 3-manifolds}, 
Topology Appl. 49 (1993), 115--127.


\bibitem{My attach}
R. Myers, {\it Attaching boundary planes to irreducible open 3-manifolds},
Quart. J. Math. Oxford Ser. (2), to appear. 

\bibitem{My endsum}
R. Myers, {\it End sums of irreducible open 3-manifolds}, 
Oklahoma State University Mathematics Department Preprint Series (1996). 

\bibitem{My cover} 
R. Myers, {\it Contractible open 3-manifolds which non-trivially cover only 
non-compact 3-manifolds}, 
Oklahoma State University Mathematics Department Preprint Series (1996). 

\bibitem{My free}
R. Myers, {\it Contractible open 3-manifolds with free covering translation 
groups}, 
Oklahoma State University Mathematics Department Preprint Series (1996). 

\bibitem{Sc}
P. Scott, {\it Compact submanifolds of 3-manifolds}, J. London Math. 
Soc. 7 (1973), 246--250. 

\bibitem{Sc-Tu} 
P. Scott and T. Tucker, {\it Some examples of exotic non-compact 
3-manifolds}, Quart. J. Math. Oxford Ser. (2) 40 (1989), 481--499. 

\bibitem{Ti-Wr}
F. Tinsley and D. Wright, {\it Some contractible open manifolds and 
coverings ofmanifolds in dimension three}, Topology Appl., to appear.

\bibitem{Tu}
T. Tucker, {\it Non-compact 3-manifolds and the missing-boundary problem}, 
Topology 13 (1974), 267--273. 

\bibitem{Wa}
F. Waldhausen, {\it On irreducible 3-manifolds which are sufficiently large},
Ann. of Math., 87 (1968), 56--88.

\bibitem{Wh}
J. H. C. Whitehead, {\it A certain open manifold whose group is unity}, 
Quart. J. Math. 6 (1935), 268--279. 



\end{thebibliography}
\end{document}